\documentclass[12pt]{article}
\usepackage{amssymb}
\usepackage{amsmath}
\usepackage{indentfirst,latexsym}
\usepackage{mathrsfs}
\usepackage{graphicx}
\usepackage{fancyhdr}
\usepackage{caption}

\date{}
\begin{document}

\title{1234-avoiding permutations and Dyck paths}
\author{Marilena Barnabei, Flavio Bonetti, and Matteo Silimbani \\ \ \\
 Dipartimento di Matematica, Universit\`a di Bologna\\ P.zza di Porta San Donato 5, 40126 Bologna, Italy} \maketitle

\noindent {\bf Abstract.} We define a map $\nu$ between the symmetric group $S_n$ and the set of pairs of Dyck paths of semilength $n$. We show that the map $\nu$ is injective when restricted to the set of $1234$-avoiding permutations and characterize the image of this map.\newline



\section{Introduction}

\noindent We say that a permutation $\sigma\in S_n$ \emph{contains a pattern}
$\tau\in S_k$ if $\sigma$ contains a subsequence that is
order-isomorphic to $\tau$. Otherwise, we say that $\sigma$
\emph{avoids} $\tau$. Given a pattern $\tau$, denote by
$S_n(\tau)$ the set of permutations in $S_n$ avoiding $\tau$.

\noindent The sets of permutations that avoid a single pattern $\tau\in S_3$ have been completely determined in last decades.
More precisely, it has been shown \cite{ss} that, for every $\tau\in S_3$, the cardinality of the set $S_n(\tau)$ equals the $n$-th Catalan number, which is also the number of Dyck paths
of semilength $n$ (see e.g. \cite{ss}). Many bijections between $S_n(\tau)$, $\tau\in S_3$, and the set of Dyck paths of semilength $n$ have been described (see \cite{ck} for a fully detailed overview).

\noindent The case of patterns of length $4$ appears much more complicated, due both to the fact that the patterns $\tau\in S_4$ are not equidistributed on $S_n$, and the difficulty of describing bijections between $S_n(\tau)$, $\tau\in S_4$, and some set of combinatorial objects.

\noindent In this paper we study the case $\tau=1234$. An explicit formula for the cardinality of $S_n(1234)$ has been computed by I.Gessel (see \cite{b} and \cite{g}), but there is no bijection (up to our knowledge) between $S_n(1234)$ and some set of combinatorial objects.

\noindent We present a bijection between $S_n(1234)$ and a set of pairs of Dyck paths of semilength $n$. More specifically, we define a map $\nu$ from $S_n$ to the set of pairs of Dyck paths, prove that every element in the image of $\nu$ corresponds to a single element in $S_n(1234)$, and characterize the set of all pairs that belong to the image of the map $\nu$.

\section{Dyck paths}
\label{pinco}

\noindent A \emph{Dyck path} of semilength $n$ is a lattice path
starting at $(0,0)$, ending at $(2n,0)$, and never going below the
$x$-axis, consisting of up steps $U=(1,1)$ and down steps
$D=(1,-1)$. A \emph{return} of a Dyck path is a down step ending
on the $x$-axis. A Dyck path is \emph{irreducible} if it has only
one return. An \emph{irreducible component} of a Dyck path $P$ is
a maximal irreducible Dyck subpath of $P$.

\noindent A Dyck path $P$ is specified by the lengths
$a_1,\ldots,a_k$ of its ascents (namely, maximal sequences of
consecutive up steps) and by the lengths  $d_1,\ldots,d_k$
 of its descents (maximal sequences of consecutive down steps), read from left to right.
 Set $A_i=\sum_{j=1}^i a_j$ and $D_i=\sum_{j=1}^i d_j$. If $n$ is the semilength of $P$,
 we have of course $A_k=D_k=n$, hence the Dyck path $P$ is uniquely determined by the two
 sequences $A=A_1,\ldots,A_{k-1}$ and $D=D_1,\ldots,D_{k-1}$. The pair $(A,D)$ is called the
 \emph{ascent-descent code} of the Dyck path $P$.

\noindent  Obviously, a pair $(A,D)$, where $A=A_1,\ldots,A_{k-1}$ and $D=D_1,\ldots,D_{k-1}$, is the ascent-descent code of
some Dyck path of semilength $n$ if and only if
\begin{itemize}
\item $0<k\leq n-1$;
\item $1\leq A_1<A_2<\ldots<A_{k-1}\leq n-1$;
\item $1\leq D_1<D_2<\ldots<D_{k-1}\leq n-1$;
\item $A_i\geq D_i$ for every $1\leq i\ \leq k-1$.
\end{itemize}
\noindent It is easy to check that the returns of a Dyck path are
in one-to-one correspondence with the indices $1\leq i\leq k$ such
that $A_i=D_i$. Hence, a Dyck path is irreducible whenever we have
$A_i>D_i$ for every $1\leq i\leq k-1$.\\

\begin{figure}[h]
\begin{center}
\includegraphics[bb=125 626 266 695,width=.5\textwidth]{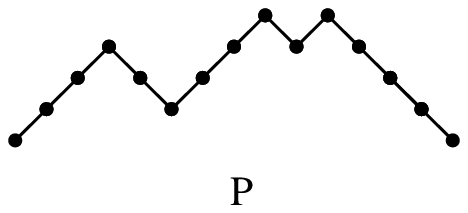}\caption{}\label{fatica}
\end{center}
\end{figure}

\noindent For example, the ascent-descent code of the Dyck path $P$ in Figure \ref{fatica} is $(A,D)$, where $A=3,6$ and $D=2,3$. Note that $A_1>D_1$ and $A_2>D_2$. In fact, $P$ is irreducible.

\begin{figure}[h]
\begin{center}
\includegraphics[bb=56 565 492 770,width=1.05\textwidth]{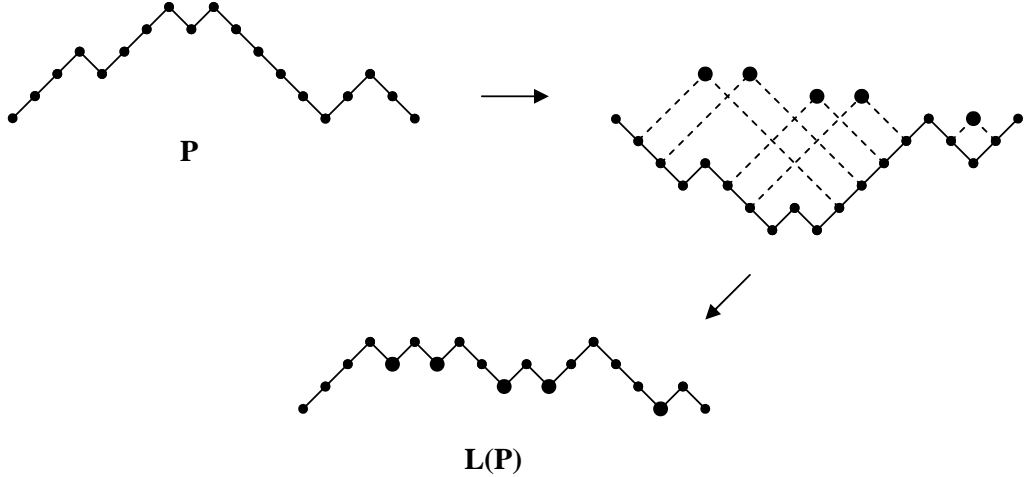}\caption{The map L.}\label{unalatra}
\end{center}
\end{figure}

\noindent We describe an involution $L$ due to Kreweras (a
description of this bijection, originally defined in \cite{kw}, can be found in \cite{c}) and discussed by Lalanne (see \cite{l1} and \cite{l2}) on
the set of Dyck paths. Given a Dyck path $P$, the path $L(P)$ can
be constructed as follows:
\begin{itemize}
\item if $P$ is the empty path $\epsilon$, than $L(P)=\epsilon$;
\item otherwise:
\begin{itemize}
\item[-] flip the Dyck path $P$ around the $x$-axis, obtaining
a path $E$; \item draw northwest (respectively northeast) lines
starting from the midpoint of each double descent (resp. ascent);
\item[-] mark the intersection between the $i$-th northwest and
$i$-th northeast line, for every $i$;
\item[-] $L(P)$ is the unique Dyck path that has valleys at
the marked points (see Figure \ref{unalatra}.
\end{itemize}
\end{itemize}

\begin{figure}[h]
\begin{center}
\includegraphics[bb=77 587 476 776,width=1.05\textwidth]{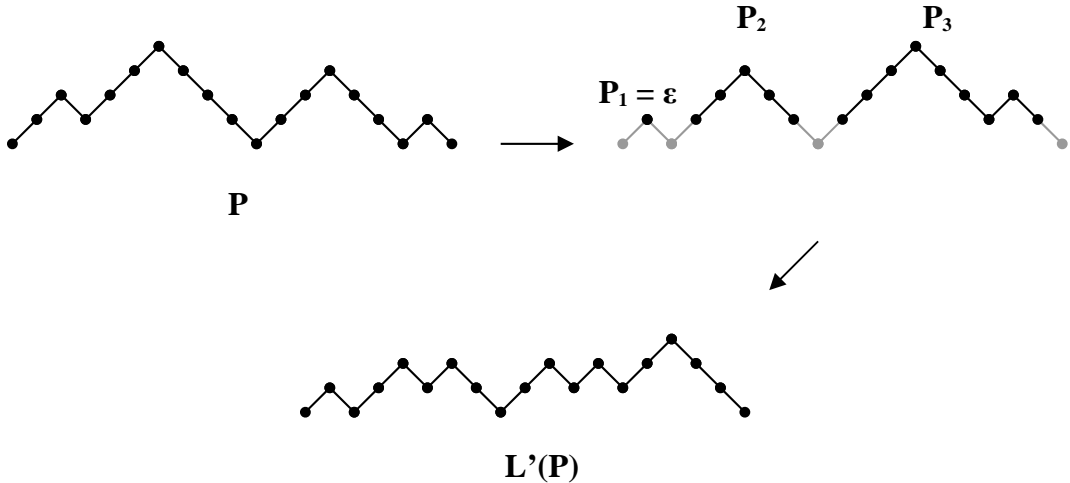}\caption{The map L'.}\label{unasola}
\end{center}
\end{figure}

\noindent We define a further involution $L'$ on the set of Dyck
paths, which is a variation of the involution $L$, as follows:
\begin{itemize}
\item if $P$ is the empty path $\epsilon$, than $L(P)=\epsilon$;
\begin{itemize}
\item[-] consider a Dyck path $P$ and flip it with respect to a vertical line; \item[-] decompose the obtained path into its irreducible
components $U\,P_i\,D$;
\item[-] replace every component $U\,P_i\,D$ with $U\,L(P_i)\,D$ in order to get $L'(P)$ (see Figure \ref{unasola}). \end{itemize}\end{itemize}

\noindent We point out that the map $L'$ appears in a slightly
modified version in the paper \cite{c}.\\

\noindent We now give a description of the map $L'$ in terms
of ascent-descent code. Obviously, it is sufficient to consider
the case of an irreducible Dyck path $P$.

\noindent Let $(A,D)$ be the ascent-descent code of an irreducible
path $P$ of semilength $n$, with $A=A_1,\ldots,A_h$ and
$D=D_1,\ldots,D_h$. Straightforward arguments show that the
ascent-descent code $(A',D')$ of $L'(P)$ can be described as
follows:
\begin{itemize}
\item set $\bar{A}_i=A_i-1$ and set $\hat{A}=[n-2]\setminus\{\bar{A}_1,\ldots,\bar{A}_h\}=\{\hat{A}_1,\ldots,\hat{A}_{n-2-h}\}$, where the $\hat{A}_i$'s are written in decreasing order.
 Then, $A'_i=n-\hat{A}_i$.
\item consider the set
$[n-2]\setminus\{D_1,\ldots,D_h\}=\{\hat{D}_1,\ldots,\hat{D}_{n-2-h}\}$,
where the $\hat{D}_i$'s are written in decreasing order. Then,
$D'_i=n-1-\hat{D}_i$.
\end{itemize}

\noindent Finally, we introduce an order relation $\leq$ on the set of Dyck paths
of the same semilength. This order relation will be defined in three steps:
\begin{itemize}
\item Consider two irreducible Dyck paths $P$
and $Q$ of semilength $n$. Let $(A,D)$ be the ascent-descent code
of $P$, with $A=A_1,\ldots,A_k$ and $D=D_1,\ldots,D_k$. We say
that $Q$ covers $P$ in the relation $\leq$ if the ascent code of $Q$ is obtained by
removing an integer $A_i$ from $A$ and the descent code of $Q$ is
obtained by removing an integer $D_j$ for $D$, with $j\geq i$.

\noindent Roughly speaking, $Q$ covers $P$ if it can be obtained
from $P$ by ``closing'' the rectangles corresponding to an
arbitrary collection of consecutive valleys of $P$;
\item the desired order relation $\leq$ on the set of
irreducible Dyck paths is the transitive closure of the above
covering relation;
\item the relation $\leq$ is extended to the set of all Dyck path of a given semilength as follows: if $P$ and $Q$ are two arbitrary Dyck paths and
$P=P_1P_2\ldots P_r$ and $Q=Q_1Q_2\ldots Q_s$ are their respective
decompositions into irreducible parts, then $P\leq Q$ whenever
$r=s$ and $P_i\leq Q_i$ for every $i$.
\end{itemize}

\begin{figure}[h]
\begin{center}
\includegraphics[bb=182 617 335 692,width=.6\textwidth]{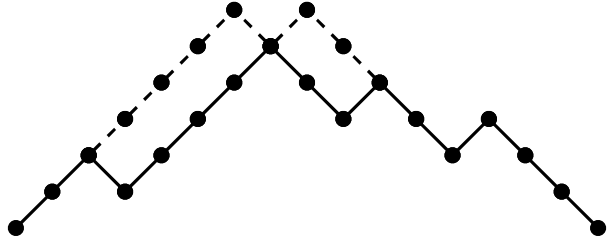}\caption{The dotted Dyck path covers the solid one.}\label{pescara}
\end{center}
\end{figure}

\section{LTR minima and RTL maxima of a permutation}

\noindent Some of
the well known bijections between $S_n(\tau)$, $\tau\in S_3$, and the set of Dyck paths of semilength $n$ (see \cite{bjs}, \cite{k}, and \cite{ss}) are based on the two notions of left-to-right minimum and right-to-left maximum of a permutation $\sigma=x_1\,x_2\,\ldots\,x_n$:

\begin{itemize}
\item
the value $x_i$ is a \emph{left-to-right minimum} (LTR minimum for
short) at position $i$ if  $x_i<x_j$ for every
$j<i$;
\item the value $x_i$ is a \emph{right-to-left maximum} (RTL
maximum) at position $i$ if $x_i>x_j$ for every
$j>i$.
\end{itemize}

For
example, the permutation
$$\sigma=5\,3\,4\,8\,2\,1\,6\,7$$  has
the LTR minima $5$, $3$, $2$, and $1$ (at positions $1$,
$2$, $5$, and $6$) and RTL
maxima $7$ and $8$ (at positions $8$ and $4$).\\

\noindent We denote by $vmin(\sigma)$ and $pmin(\sigma)$ the sets
of values and positions of the LTR minima of $\sigma$, respectively.
Analogously, $vmax(\sigma)$ and $pmax(\sigma)$ denote the sets of
values and positions of the RTL maxima of $\sigma$.\\

\noindent Recall that the \emph{reverse-complement} of a permutation
$\sigma\in S_n$ is the permutation defined by
$$\sigma^{rc}(i)=n+1-\sigma(n+1-i).$$
For example, consider the permutation
$$\sigma=2\, 4\, 7\, 3\, 1\, 8\, 9\, 5\, 6.$$ Then:
$$\sigma^{rc}=4\, 5\, 1\, 2\, 9\, 7\, 3\, 6\, 8.$$

\noindent Note that the sets $S_n(123)$ and $S_n(1234)$ are closed
under reverse-complement, namely, $\sigma\in S_n(123)$
(respectively, $\sigma\in S_n(1234)$) if and only if
$\sigma^{rc}\in S_n(123)$ (resp. $\sigma^{rc}\in S_n(1234)$).\\

\noindent The first assertion in the next proposition goes back to the seminal paper \cite{ss},
while the second one is an immediate consequence of the straightforward fact that
$x$ is a LTR minimum of a permutation $\sigma$ at position $i$ if and only if $n+1-x$ is
RTL maximum of $\sigma^{rc}$ at position $n+1-i$:

\newtheorem{iniziale}{Theorem}
\begin{iniziale}\label{possiamo} A permutation $\sigma\in S_n(123)$ is completely determined by the two sets
$vmin(\sigma)$ and $pmin(\sigma)$ of values and positions of its
left-to-right minima. A permutation in $S_n(123)$ is
completely determined, as well, by the two sets $vmax(\sigma)$ and
$pmax(\sigma)$ of values and positions of its right-to-left
maxima.\end{iniziale}
\begin{flushright}
\vspace{-.4cm}$\diamond$
\end{flushright}

\noindent Also $1234$-avoiding permutations can be characterized
in terms of LTR minima and RTL maxima.

\noindent This characterization can be found in \cite{b} and is
based on an equivalence relation on $S_n$ defined as follows:
$\sigma\equiv\sigma'$ $\iff$ $\sigma$ and $\sigma'$ share the
values and the positions of LTR minima and RTL maxima.\newline

\noindent For example,
$$1\,2\,3\,4\equiv 1\,3\,2\,4.$$
Straightforward arguments lead to the following result stated in
\cite{b}:

\newtheorem{bona}[iniziale]{Theorem}
\begin{bona} \label{le} Every  equivalence class of the relation $\equiv$ contains exactly one
$1234$-avoiding permutation. In this permutation, the values that
are neither LTR minima nor RTL maxima appear in decreasing order.
\end{bona}
\begin{flushright}
\vspace{-.4cm}$\diamond$
\end{flushright}

\section{The maps $\lambda$ and $\mu$}

\noindent We define two maps $\lambda$ and $\mu$ between $S_n$ and
the set $\mathcal{D}_n$ of Dyck paths of semilength $n$. Given a
permutation $\sigma\in S_n$, the path $\lambda(\sigma)$ is
contructed as follows:

\begin{itemize}
\item decompose $\sigma$ as
$\sigma=m_1\,w_1\,m_2\,w_2\,\ldots\,m_k\,w_k$, where
$m_1,m_2,\ldots,m_k$ are the left-to-right minima in $\sigma$ and
$w_1,w_2,\ldots,w_k$ are (possibly empty) words;
\item  set $m_0=n+1$;
\item read the permutation from left to right and translate any
LTR minimum $m_i$ ($i>0$) into $m_{i-1}-m_{i}$ up steps and any
subword $w_i$ into $l_i+1$ down steps, where $l_i$ denotes the
number of elements in $w_i$.
\end{itemize} The statement of Theorem \ref{possiamo} implies that the map $\lambda$ is a bijection when restricted to $S_n(123)$.

\noindent Note that the ascent-descent code $(A,D)$ of the path
$\lambda(\sigma)$ is obtained as follows: \begin{itemize}
\item $A=n+1-m_1,n+1-m_2,\ldots,n+1-m_{k-1}$;
\item $D=p_2-1,p_3-1,\ldots,p_k-1$, where $p_i$ is the position of
$m_i$.
\end{itemize}

\noindent We
define a further map $\mu:S_n\to\mathcal{D}_n$:

\begin{itemize}
\item decompose $\sigma$ as
$\sigma=u_h\,M_h\,u_{h-1}\,M_{h-1}\,\ldots\,u_1\,M_1$, where
$M_1,M_2,\ldots,M_h$ are the right-to-left maxima in $\sigma$ and
$u_1,u_2,\ldots,u_k$ are (possibly empty) words;\item set $M_0=0$; \item associate with
$M_i$ ($i>0$) the steps $U^{m_i-m_{i-1}}D$
\item associate with each entry in $u_i$ a $D$ step.
\end{itemize}
 Also in this case, the map $\mu$ is a bijection when restricted to
$S_n(123)$.\\

\noindent The ascent-descent code $(A^*,D^*)$ of
the path $\mu(\sigma)$ is obtained as follows: \begin{itemize}
\item $A^*=M_1,M_2,\ldots,M_{h-1}$;
\item $D^*=n-P_2,n-P_3,\ldots,n-P_h$, where $P_i$ is the position of
$M_i$.
\end{itemize}

\noindent In Figure \ref{pescara} the two paths $\lambda(\sigma)$ and $\mu(\sigma)$ corresponding to $\sigma=6\,2\,3\,1\,7\,5\,4$ are shown.

\begin{figure}[h]
\begin{center}
\includegraphics[bb=71 688 377 764,width=.9\textwidth]{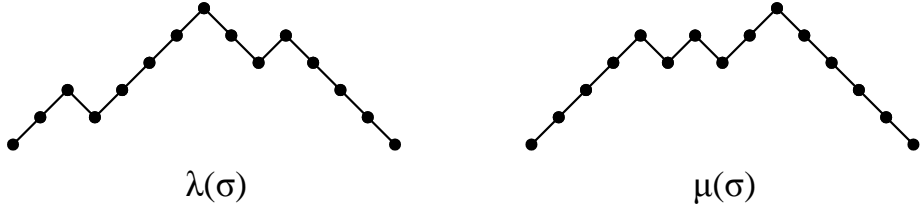}\caption{The Dyck paths corresponding to $\sigma=6\,2\,3\,1\,7\,5\,4$.}\label{pescara}
\end{center}
\end{figure}

\noindent We can now define a map
$\nu:S_n\to\mathcal{D}_n\times\mathcal{D}_n$, setting
$$\nu(\sigma)=(\lambda(\sigma),\mu(\sigma)).$$

\noindent The statement of Theorem \ref{le} implies that the map $\nu$ is injective when
restricted to $S_n(1234)$.\newline

\noindent Note that the map $\nu$ behaves properly with respect to
the reverse-complement and the inversion operators:

\newtheorem{notache}[iniziale]{Proposition}
\begin{notache} Let $\sigma$ be a permutation in $S_n$. We
have:\begin{itemize} \item $\nu(\sigma)=(L,R)\iff
\nu(\sigma^{rc})=(R,L),$ hence, the permutation $\sigma$ is rc-invariant if and
only if $L=R$.
\item $\nu(\sigma)=(L,R)\iff\nu(\sigma^{-1})=(rev(L),rev(R))$, where $rev(P)$ is the path
obtained by flipping $P$ with respect to a vertical line. Hence, the permutation
 $\sigma$ is an involution if and only if both $L$
and $R$ are symmetric with respect to a vertical line.\newline
\end{itemize}
\end{notache}
\begin{flushright}
\vspace{-.4cm}$\diamond$
\end{flushright}

\begin{figure}[h]
\begin{center}
\includegraphics[bb=71 688 377 764,width=.9\textwidth]{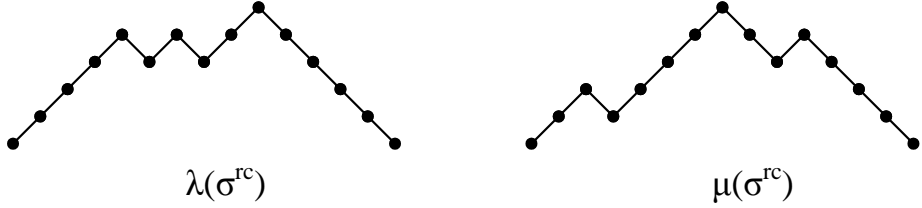}\caption{The Dyck paths corresponding to $\sigma^{rc}=4\,3\,1\,7\,5\,6\,2$.}\label{cocco}
\end{center}
\end{figure}

\begin{figure}[h]
\begin{center}
\includegraphics[bb=71 688 377 764,width=.9\textwidth]{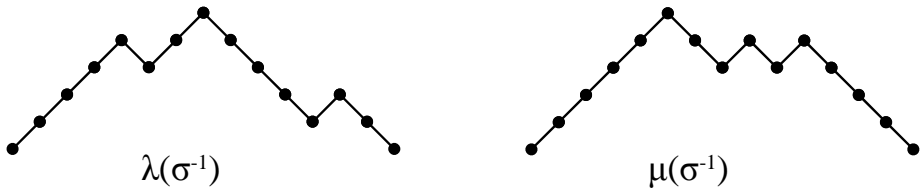}\caption{The Dyck paths corresponding to $\sigma^{-1}=4\,2\,3\,7\,6\,1\,5$.}\label{marilena}
\end{center}
\end{figure}

\noindent For example, consider
$\sigma=6\,2\,3\,1\,7\,5\,4$. The two paths associated with $\sigma$ are shown in Figure \ref{pescara}.
The permutation $\sigma^{rc}=4\,3\,1\,7\,5\,6\,2$
is associated with the two paths in Figure \ref{cocco}, while the permutation $\sigma^{-1}=4\,2\,3\,7\,6\,1\,5$ corresponds to the two paths in Figure \ref{marilena}.\\

\noindent Moreover, the map $\nu$ has the following further
property that will be crucial in the proof of our main result.

\noindent Recall that a permutation $\sigma\in S_n$ is said \emph{right-connected} if
it does not have a suffix $\sigma'$ of length $k<n$, that is a
permutation of the symbols $1,2,\ldots,k$.

\noindent For example, the permutation
$$\tau=6\,1\,2\,7\,5\,3\,4\,8$$ is right-connected, while $$\sigma=8\,6\,4\,5\,7\,2\,1\,3$$
is not right-connected.

\noindent According to this definition, we can split every
permutation into right-connected components:
$$\sigma={\bf 8}\,6\,4\,5\,7\,{\bf 2\,1\,3}.$$

\noindent Note that, if a permutation $\sigma$ is not right-connected, $\sigma$ is the juxtaposition of a
permutation $\sigma''$ of the set $\{t+1,\ldots,n\}$ and the permutation
$\sigma'$ of the set $\{1,\ldots,t\}$.

\newtheorem{componenti}[iniziale]{Proposition}
\begin{componenti}\label{scusa}
Let $\sigma$ be a non right-connected permutation in $S_n$, with $\sigma=\sigma_1\sigma_2$, where $\sigma_1$ is a permutation of the set $\{t+1,\ldots,n\}$ and
$\sigma_2$ is a permutation of set of the set $\{1,\ldots,t\}$. Then:
$$\lambda(\sigma)=P_1P_2\qquad \mu(\sigma)=Q_1Q_2,$$
with $P_i=\lambda(\sigma_i)$ and $Q_i=\mu(\sigma_i)$, $i=1,2$.
\end{componenti}
\begin{flushright}
\vspace{-.4cm}$\diamond$
\end{flushright}

\noindent The order relation on Dyck paths defined in Section \ref{pinco} can be
exploited to define two order relations on the set $S_n$ as
follows: \begin{itemize} \item $\sigma\leq_{\lambda} \tau$ if and
only if $\lambda(\sigma)\leq\lambda(\tau)$; \item
$\sigma\leq_{\mu}\tau$ if and only if $\mu(\sigma)\leq
\mu(\tau)$.\end{itemize} These order relations can be
intrinsically described as follows:

\newtheorem{hallo}[iniziale]{Proposition}
\begin{hallo} \label{ancorano} Let $\sigma,\tau\in S_n$. We have $\sigma\leq_{\lambda} \tau$ whenever:
\begin{itemize}
\item $vmin(\tau)\subseteq vmin(\sigma)$;
\item $pmin(\tau)\subseteq pmin(\sigma)$;
\item setting:\\ $vmin(\sigma)=\{m_1,\ldots,m_h\}$ (written in decreasing order),\\
$vmin(\sigma)\setminus
vmin(\tau)=\{m_{i_1},m_{i_2},\ldots,m_{i_r}\}$ (in decreasing
order),\\ $pmin(\sigma)\setminus
pmin(\tau)=\{p_{j_1},p_{j_2},\ldots,p_{j_r}\}$ (in increasing
order),\\ then $i_k<j_k$ for every $k$.
\end{itemize}
Similarly, $\sigma\leq_{\mu} \tau$ whenever:
\begin{itemize}
\item $vmax(\tau)\subseteq vmax(\sigma)$;
\item $pmax(\tau)\subseteq pmax(\sigma)$;
\item setting:\\ $vmax(\sigma)=\{M_1,\ldots,M_t\}$ (written in increasing order),\\ $vmax(\sigma)\setminus vmax(\tau)=
\{M_{i_1},M_{i_2},\ldots,M_{i_q}\}$ (in increasing order),\\
$pmax(\sigma)\setminus pmax(\tau)=
\{P_{j_1},P_{j_2},\ldots,P_{j_q}\}$ (in decreasing order),\\ then
$i_k<j_k$ for every $k$.
\end{itemize}
\end{hallo}
\begin{flushright}
\vspace{-.4cm}$\diamond$
\end{flushright}

\noindent For example, consider the permutation
$$\sigma=6\,8\,7\,3\,2\,5\,9\,1\,4.$$
We have $vmin(\sigma)=\{6,3,2,1\}$, $pmin(\sigma)=\{1,4,5,8\}$,
$vmax(\sigma)=\{4,9\}$, and $pmax(\sigma)=\{9,7\}$. The
permutation $$\tau=3\,4\,9\,2\,6\,8\,7\,1\,5$$ is such that
$vmin(\tau)=\{3,2,1\}$ and $pmin(\tau)=\{1,4,8\}$, hence,
$\sigma\leq_{\lambda} \tau$. Moreover, the permutation
$$\rho=2\,7\,1\,3\,4\,6\,5\,8\,9$$ is such that
$vmax(\rho)=\{9\}$ and $pmax(\rho)=\{9\}$, hence,
$\sigma\leq_{\mu} \rho$.

\section{Main results}

\noindent We say that a pair of Dyck paths $(P,Q)$ is
\emph{admissible} if there exists a permutation $\alpha$ such that
$P=\lambda(\alpha)$ and $Q=\mu(\alpha)$. Needless to say, the set
of admissible pairs is in bijection with the set of
$1234$-avoiding permutations.

\noindent We want to show that the operator $L'$ on Dyck paths
allows us to characterize the set of admissible pairs. We begin
with a preliminary result concerning the pairs of Dyck paths
corresponding to $123$-avoiding permutations:

\newtheorem{relazione}[iniziale]{Theorem}
\begin{relazione} For every $\sigma\in
S_n(123)$, we have:$$\mu(\sigma)=L'(\lambda(\sigma)).$$
\end{relazione}

\noindent\emph{Proof} Proposition \ref{scusa}, together with the
definition of the map $L'$, allows us to restrict our attention to
the right-connected case.

\noindent Recall (see \cite{ss}) that a permutation $\sigma$
avoids $123$ if and only if the set $vmin(\sigma)\cup
vmax(\sigma)=[n]$. It is simple to check that, if $\sigma$ is
right-connected, the sets of LTR minima and RTL maxima are
disjoint.

\noindent Consider now a permutation $\sigma$ with LTR minima
$m_1,\ldots,m_{k-1},m_k=1$ and RTL maxima $M_1,\ldots,M_{h-1},M_h=n$. Denote by
$(A,D)$ the ascent-descent code of the path $P=\lambda(\sigma)$
and by $(A^*,D^*)$ the ascent-descent code of the path
$\mu(\sigma)$.

\noindent As noted before, the ascent code $A'$ of $L'(P)$ is
obtained by computing the integers $\bar{A}_i=A_i-1$ and then
considering the set
$\hat{A}=[n-2]\setminus\{\bar{A}_1,\ldots,\bar{A}_{k-1}\}$, which
can be written as
$$\hat{A}=\{n-(n-1),n-(n-2),\ldots,n-2\}\setminus\{n-m_1,\ldots,n-m_{k-1}\}.$$
Since
$\{m_1,\ldots,m_{k-1}\}\cup\{M_1,\ldots,M_{h-1}\}=\{2,3,\ldots,n-1\}$,
we have
$$\hat{A}=\{n-M_1,\ldots,n-M_{h-1}\}.$$
Hence, $A'=A^*$.

\noindent Similarly, the descent code $D'$ of $L'(P)$ is obtained
by considering the set
$$\hat{D}=[n-2]\setminus\{D_1,\ldots,D_{k-1}\}=[n-2]\setminus\{p_2-1,\ldots,p_k-1\}.$$
Since
$\{p_1,\ldots,p_{k-1}\}\cup\{P_1,\ldots,P_{h-1}\}=\{2,3,\ldots,n-1\}$,
we have $$\hat{D}=\{P_2-1,\ldots,P_{h-1}-1\}.$$ Hence, $D'=D^*$.
\begin{flushright}
\vspace{-.4cm}$\diamond$
\end{flushright}

\noindent For example, the $123$-avoiding permutation $\sigma=8\,5\,9\,7\,6\,2\,4\,3\,1$ corresponds to the pair of Dyck paths $(P,L'(P))$ in Figure \ref{unasola}.\\

\noindent We are now in position to state our main result:
\newtheorem{finale}[iniziale]{Theorem}
\begin{finale}
A pair $(P,Q)$ is admissible if and only if $P\geq L'(Q)$ and
$Q\geq L'(P)$.
\end{finale}

\noindent\emph{Proof} Consider a permutation $\sigma\in S_n(1234)$
and let $\sigma'$ be the unique permutation in $S_n(123)$ with the
same LTR minima as $\sigma$, at the same positions. Obviously,
$\sigma'\leq_{\mu}\sigma$, since in $\sigma'$ every element that
is not a LTR minimum is a RTL maximum (see Proposition
\ref{ancorano}). Recalling that
$\mu(\sigma')=L'(\lambda(\sigma))=L'(P)$, we get the first
inequality. The other inequality follows from the fact that the
pair $(P,Q)$ is admissible whenever the pair $(Q,P)$ is
admissible.\newline

\noindent Consider now a pair of Dyck paths $(P,Q)$ such that
$P\geq L'(Q)$ and $Q\geq L'(P)$. Proposition
\ref{scusa} allows us to restrict to the case $P,Q$ irreducible.
Denote by $\sigma$ and $\tau$ the permutations in $S_n(123)$
corresponding via $\nu$ to the pairs $(P,L'(P))$ and $(L'(Q),Q)$,
respectively. Since $P\geq L'(Q)$ and $Q\geq L'(P)$, we
have $\tau\leq_{\lambda}\sigma$ and $\sigma\leq_{\mu}\tau$.\\

\noindent We define a permutation $\alpha\in S_n$ as follows:
\begin{itemize}
\item $\alpha(x)=\sigma(x)$ if $x\in pmin(\sigma)$;
\item $\alpha(x)=\tau(x)$ if $x\in pmax(\tau)$;
\item if $x\notin pmin(\sigma)\cup pmax(\tau)$,
we have $x\in pmax(\sigma)\setminus pmax(\tau)=pmin(\tau)\setminus
pmin(\sigma)=\{p_{j_1},\ldots,p_{j_r}\}$, written in increasing
order. Set $$\alpha(p_{j_k})=m_{i_k},$$ where
$m_{i_1},m_{i_2},\ldots,m_{i_r}$ are the elements in
$vmin(\tau)\setminus vmin(\sigma)=vmax(\sigma)\setminus
vmax(\tau)$, written in decreasing order.
\end{itemize}
The permutation $\alpha$ is obtained as the concatenation of three
decreasing sequences. Hence, $\alpha$ avoids $1234$. We have to
prove that $vmin(\alpha)=vmin(\sigma)$ and
$vmax(\alpha)=vmax(\tau)$.

\noindent It is immediate that $vmin(\sigma)\subseteq
vmin(\alpha)$. In order to prove that $vmin(\sigma)= vmin(\alpha)$
it remains to show that the values
$m_{i_1},m_{i_2},\ldots,m_{i_r}$ are not LTR minima of $\alpha$.

\noindent In fact, for every $k$, consider
$\alpha(p_{j_k})=m_{i_k}=\tau(p_{i_k})$. Consider the sets
$A=\{p_1,p_2,\ldots,p_{i_k}\}$, $B=\{m_1,m_2,\ldots,m_{i_k}\}$, and
their subsets $A'=\{p_{i_1},p_{i_2},\ldots,p_{i_k}\}$ and
$B'=\{m_{i_1},m_{i_2},\ldots,m_{i_k}\}$. The $k$ elements in $B'$
do not belong to $vmin(\sigma)$ (and hence, the $i_k-k$ elements
in $B\setminus B'$ are the largest elements in $vmin(\sigma)$).
Proposition \ref{ancorano} ensures that each of them occupies in
$\alpha$ a position that is strictly greater than the position
occupied in $\tau$. This implies that $p_{j_k}<p_{i_k}$ and that
at most $k-1$ elements in $B'$ occupy in $\tau$ a position that belongs to $A$.
Hence, in $\alpha$, at least $i_k-k+1$ positions in $A$ are
occupied by entries belonging to $vmin(\sigma)$. This implies that
there is in $\alpha$ a position preceding $p_{j_k}$ occupied by a
value less than $m_{i_k}$. Hence, $m_{i_k}$ is not a LTR minimum
of $\alpha$.

\noindent Analogous arguments can be used to prove that
$vmax(\alpha)=vmax(\tau)$. Hence, $\nu(\alpha)=(P,Q)$, as desired.
\begin{flushright}
\vspace{-.4cm}$\diamond$
\end{flushright}

\begin{figure}[h]
\begin{center}
\includegraphics[bb=107 590 502 680,width=.9\textwidth]{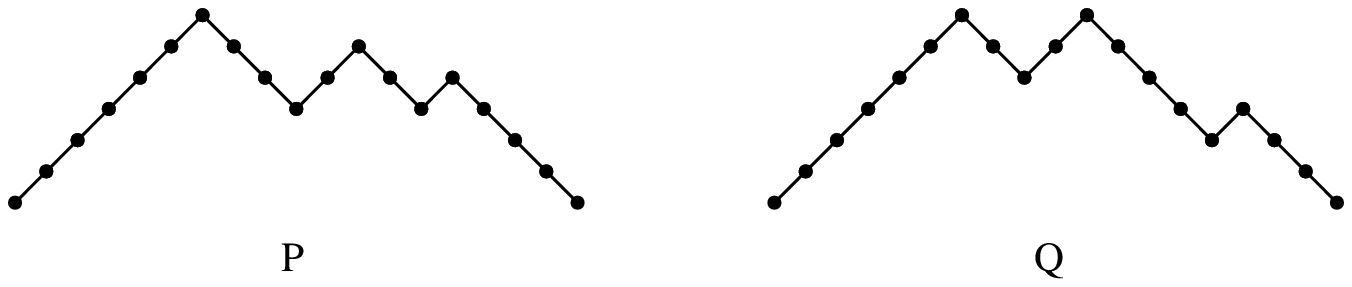}\caption{}\label{tizcaj}
\end{center}
\end{figure}

\noindent For example, consider the pair of Dyck paths in Figure \ref{tizcaj}.\\

\noindent It can
be checked that $P\geq L'(Q)$ and $Q\geq L'(P)$. The permutations
$\sigma=\nu^{-1}((P,L'(P)))$ and $\tau=\nu^{-1}((L'(Q),Q))$ are as
follows:
$$\sigma=4\,9\,8\,2\,7\,1\,6\,5\,3\qquad\tau=7\,5\,9\,4\,3\,2\,8\,1\,6.$$
We have $vmin(\sigma)=\{4,2,1\}$, $pmin(\sigma)=\{1,4,6\}$,
$vmin(\tau)=\{7,5,4,3,2,1\}$, $pmin(\tau)=\{1,2,4,5,6,8\}$,
$vmax(\sigma)=\{3,5,6,7,8,9\}$, $pmax(\sigma)=\{9,8,7,5,3,2\}$,
$vmax(\tau)=\{6,8,9\}$, and $pmax(\tau)=\{9,7,3\}$.

\noindent The
permutation $\alpha=\nu^{-1}((P,Q))$ is
$$\alpha=4\,7\,9\,2\,5\,1\,8\,3\,6.$$
As expected, $vmin(\alpha)=vmin(\sigma)$,
$pmin(\alpha)=pmin(\sigma)$, $vamx(\alpha)=vmax(\tau)$, and
$pmax(\alpha)=pmax(\tau)$.

\end{document}